\documentclass[12pt]{article}%
\usepackage{graphicx}
\usepackage[intlimits]{amsmath}
\usepackage{latexsym}
\usepackage{amsfonts}
\usepackage{amssymb}%
 \makeatletter
\newfont{\footsc}{cmcsc10 at 8truept}
\newfont{\footbf}{cmbx10 at 8truept}
\newfont{\footrm}{cmr10 at 10truept}
\newtheorem{theorem}{Theorem}

\newtheorem{lemma}[theorem]{Lemma}

\newenvironment{proof}[1][Proof]{\noindent{\textbf {#1}  }}  {\hfill$\Box$\bigskip}

\begin{document}

\title{Graphs with many copies of a given subgraph}
\author{Vladimir Nikiforov\\{\small Department of Mathematical Sciences, University of Memphis, Memphis TN
38152}}
\maketitle

\begin{abstract}
Let $c>0,$ and $H$ be a fixed graph of order $r.$ Every graph on $n$ vertices
containing at least $cn^{r}$ copies of $H$ contains a \textquotedblleft
blow-up\textquotedblright\ of $H$ with $r-1$ vertex classes of size
$\left\lfloor c^{r^{2}}\ln n\right\rfloor $ and one vertex class of size
greater than $n^{1-c^{r-1}}.$ A similar result holds for induced copies of
$H.$\medskip


\end{abstract}

\subsection*{Main results}

This note is part of an ongoing project aiming to renovate some classical
results in extremal graph theory, see, e.g., \cite{BoNi04} and
\cite{Nik07,Nik07c}.

Suppose that a graph $G$ of order $n$ contains $cn^{r}$ copies of a given
subgraph $H$ on $r$ vertices. How large \textquotedblleft
blow-up\textquotedblright\ of $H$ must $G$ contain? When $H$ is an $r$-clique,
this question was answered in \cite{Nik07}: $G$ contains a complete
$r$-partite graph with $r-1$ parts of size $\left\lfloor c^{r}\ln
n\right\rfloor $ and one part larger than $n^{1-c^{r-1}}$.

The aim of this note is to answer this question for any subgraph $H$.

We first define precisely a \textquotedblleft blow-up\textquotedblright\ of a
graph: given a graph $H$ of order $r$ and positive integers $x_{1}%
,\ldots,x_{r}$, we write $H\left(  x_{1},\ldots,x_{r}\right)  $ for the graph
obtained by replacing each vertex $u\in V\left(  H\right)  $ with a set
$V_{u}$ of size $x_{u}$ and each edge $uv\in E\left(  H\right)  $ with a
complete bipartite graph with vertex classes $V_{u}$ and $V_{v}.$

\begin{theorem}
\label{th1}Let $r\geq2,$ $0<c\leq1/2,$ $H$ be a graph of order $r,$ and $G$ be
a graph of order $n.$ If $G$ contains more than $cn^{r}$ copies of $H,$ then
$H\left(  s,\ldots s,t\right)  \subset G,$ where $s=\left\lfloor c^{r^{2}}\ln
n\right\rfloor $ and $t>n^{1-c^{r-1}}.$
\end{theorem}

To state a similar theorem for induced subgraphs, we introduce a new concept:
we say that a graph $X$ \emph{is of type }$H\left(  x_{1},\ldots,x_{r}\right)
,$ if $X$ is obtained from $H\left(  x_{1},\ldots,x_{r}\right)  $ by adding
some edges within the sets $V_{u},$ $u\in V\left(  H\right)  .$

\begin{theorem}
\label{th2}Let $r\geq2,$ $0<c\leq1/2,$ $H$ be a graph of order $r,$ and $G$ be
a graph of order $n.$ If $G$ contains more than $cn^{r}$ induced copies of
$H,$ then $G$ contains an induced subgraph of type $H\left(  s,\ldots
s,t\right)  ,$ where $s=\left\lfloor c^{r^{2}}\ln n\right\rfloor $ and
$t>n^{1-c^{r-1}}.$
\end{theorem}

The proofs of Theorems \ref{th1} and \ref{th2} are almost identical; we shall
present only the proof of Theorem \ref{th2}, for it needs more care.

Our notation follows \cite{Bol98}. Thus $V\left(  G\right)  $ and $E\left(
G\right)  $ denote the vertex and edge sets of a graph $G$ and $e\left(
G\right)  =\left\vert E\left(  G\right)  \right\vert .$ The subgraph induced
by $X\subset V\left(  G\right)  $ is denoted by $G\left[  X\right]  .$\bigskip

\textbf{Specific notation}\bigskip

Suppose $G$ and $H$ are graphs, and $X$ is an induced subgraph of $H.$

We write $H\left(  G\right)  $ for the set of injections $h:H\rightarrow G,$
such that $\left\{  u,v\right\}  \in E\left(  H\right)  $ if and only if
$\left\{  h\left(  u\right)  ,h\left(  v\right)  \right\}  \in E\left(
G\right)  .$

We say that $P\in H\left(  G\right)  $ \emph{extends} $R\in X\left(  G\right)
,$ if $R=P|V\left(  X\right)  .$

Suppose $M\subset H\left(  G\right)  .$

We let
\[
X\left(  M\right)  =\left\{  R:\left(  R\in X\left(  G\right)  \right)  \text{
}\&\text{ }\left(  \text{there exists }P\in M\text{ extending }R\right)
\right\}  .
\]

For every $R\in X\left(  M\right)  ,$ we let
\[
d_{M}\left(  R\right)  =\left\vert \left\{  P:\left(  P\in M\right)  \text{
}\&\text{ }\left(  P\text{ extends }R\right)  \right\}  \right\vert .
\]

Suppose $Y$ is a subgraph of $G$ of type $H\left(  s_{1},\ldots,s_{r}\right)
$ and $s=\min\left\{  s_{1},\ldots,s_{r}\right\}  .$

We say that $M$ \emph{covers} $Y$ if:

(a) for every edge $ij$ going across vertex classes of $Y,$ there exists $h\in
K_{2}\left(  M\right)  $ mapping some edge of $H$ onto $ij;$

(b) there exists $h_{1},\ldots,h_{s}\in M,$ such that $h_{i}\left(  H\right)
\cap h_{j}\left(  H\right)  =\varnothing$ for $i\neq j,$ and for all
$i\in\left[  s\right]  ,$ $h_{i}\left(  H\right)  $ intersects all vertex
classes of $Y.$\bigskip

We deduce Theorem \ref{th2} from the following technical statement.

\begin{theorem}
\label{genZ}Let $r\geq2,$ $0<c\leq1/2,$ $H$ be a graph of order $r,$ and $G$
be a graph of order $n.$ If $M\subset H\left(  G\right)  $ and $\left\vert
M\right\vert \geq cn^{r},$ then $M$ covers an induced subgraph of type
$H\left(  s,\ldots s,t\right)  $ with $s=\left\lfloor c^{r^{2}}\ln
n\right\rfloor $ and $t>n^{1-c^{r-1}}.$
\end{theorem}

The proof of Theorem \ref{genZ} is based on the following routine lemma.

\begin{lemma}
\label{le1} Let $F$ be a bipartite graph with parts $A$ and $B.$ Let
$\left\vert A\right\vert =m,$ $\left\vert B\right\vert =n,$ $r\geq2,$
$0<c\leq1/2,$ and $s=\left\lfloor c^{r^{2}}\ln n\right\rfloor .$ If
$s\leq\left(  c/2^{r}\right)  m+1$ and $e\left(  F\right)  \geq\left(
c/2^{r-1}\right)  mn,$ then $F$ contains a $K_{2}\left(  s,t\right)  $ with
parts $S\subset A$ and $T\subset B$ such that $\left\vert S\right\vert =s$ and
$\left\vert T\right\vert =t>n^{1-c^{r-1}}$.
\end{lemma}

\begin{proof}
Let
\[
t=\max\left\{  x:\text{there exists }K_{2}\left(  s,x\right)  \subset F\text{
with part of size }s\text{ in }A\right\}  .
\]

For any $X\subset A,$ write $d\left(  X\right)  $ for the number of vertices
joined to all vertices of $X.$ By definition, $d\left(  X\right)  \leq t$ for
each $X\subset A$ with $\left\vert X\right\vert =s;$ hence,%
\begin{equation}
t\binom{m}{s}\geq%
{\textstyle\sum\limits_{X\subset A,\left\vert X\right\vert =s}}
d\left(  X\right)  =%
{\textstyle\sum\limits_{u\in B}}
\binom{d\left(  u\right)  }{s}. \label{in1}%
\end{equation}
Following \cite{Lov79}, p. 398, set
\[
f\left(  x\right)  =\left\{
\begin{array}
[c]{cc}%
\binom{x}{s} & \text{if }x\geq s-1\\
0 & \text{if }x<s-1,
\end{array}
\right.
\]
and note that $f\left(  x\right)  $ is a convex function. Therefore,%
\[%
{\textstyle\sum\limits_{u\in B}}
\binom{d\left(  u\right)  }{s}=%
{\textstyle\sum\limits_{u\in B}}
f\left(  d\left(  u\right)  \right)  \geq nf\left(  \frac{1}{n}%
{\textstyle\sum\limits_{u\in B}}
d\left(  u\right)  \right)  =n\binom{e\left(  F\right)  /n}{s}\geq
n\binom{cm/2^{r-1}}{s}.
\]
Combining this inequality with (\ref{in1}), and rearranging, we find that%
\begin{align*}
t  &  \geq n\frac{\left(  cm/2^{r-1}\right)  \left(  cm/2^{r-1}-1\right)
\cdots\left(  cm/2^{r-1}-s+1\right)  }{m\left(  m-1\right)  \cdots\left(
m-s+1\right)  }>n\left(  \frac{cm/2^{r-1}-s+1}{m}\right)  ^{s}\\
&  \geq n\left(  \frac{c}{2^{r}}\right)  ^{s}\geq n\left(  e^{\ln\left(
c/2^{r}\right)  }\right)  ^{c^{r^{2}}\ln n}=n^{1+c^{r^{2}}\ln\left(
c/2^{r}\right)  }=n^{1+2^{r}c^{r^{2}-1}\left(  c/2^{r}\right)  \ln\left(
c/2^{r}\right)  }.
\end{align*}
Since $c/2^{r}\leq1/8<1/e$ and $x\ln x$ is decreasing for $0<x<1/e,$ we see
that%
\[
t>n^{1+2^{r}c^{r^{2}-1}\left(  c/2^{r}\right)  \ln\left(  c/2^{r}\right)
}\geq n^{1-\left(  r+1\right)  c^{r^{2}-1}\left(  1/2\right)  \ln
2}>n^{1-c^{r-1}\left(  r+1\right)  2^{-r^{2}+r}\left(  1/2\right)  \ln2}.
\]
Now, $t>n^{1-c^{r-1}}$ follows, in view of
\[
\frac{2^{r^{2}-r}}{r+1}\geq\frac{1}{2\ln2},
\]
completing the proof.
\end{proof}

\bigskip

\begin{proof}
[\textbf{Proof of Theorem \ref{genZ}}]Let $M\subset H\left(  G\right)  $
satisfy $\left\vert M\right\vert \geq cn^{r}.$ We shall use induction on $r$
to prove that $M$ covers an induced subgraph of type $H\left(  s,\ldots
s,t\right)  $ with $s=\left\lfloor c^{r^{2}}\ln n\right\rfloor $ and
$t>n^{1-c^{r-1}}.$

Assume $r=2$ and let $A$ and $B$ be two disjoint copies of $V\left(  G\right)
.$ We can suppose that $H=K_{2},$ as otherwise we apply the subsequent
argument to the complement of $G.$

Define a bipartite graph $F$ with parts $A$ and $B,$ joining $u\in A$ to $v\in
B$ if $uv\in M.$ Set $s=\left\lfloor c^{4}\ln n\right\rfloor $ and note that
$s\leq\left(  c/4\right)  n+1.$ Since $e\left(  F\right)  =\left\vert
M\right\vert \geq cn^{2}>\left(  c/2\right)  n^{2},$ Lemma \ref{le1} implies
that $F$ contains a $K_{2}\left(  s,t\right)  $ with $t>n^{1-c}.$ Hence $M$
covers an induced graph of type $K_{2}\left(  s,t\right)  ,$ proving the
assertion for $r=2.$ Assume the assertion true for $2\leq r^{\prime}<r.$

Let $V\left(  H\right)  =\left\{  v_{1},\ldots,v_{r}\right\}  $ and
$H^{\prime}=H\left[  \left\{  v_{1},\ldots,v_{r-1}\right\}  \right]  .$

We first show that there exists $L\subset M$ with $\left\vert L\right\vert
>\left(  c/2\right)  n^{r}$ such that $d_{L}\left(  R\right)  >\left(
c/2\right)  n$ for all $R\in H^{\prime}\left(  L\right)  .$ Indeed, set $L=M$
and apply the following procedure.

\textbf{While }\emph{there exists an }$R\in H^{\prime}\left(  L\right)
$\emph{ with }$d_{L}\left(  R\right)  \leq\left(  c/2\right)  n$ \textbf{do}

\qquad\emph{Remove from }$L$\emph{ all members extending }$R.$

When this procedure stops, we have $d_{L}\left(  R\right)  >\left(
c/2\right)  n$ for all $R\in H^{\prime}\left(  L\right)  ,$ and also
\[
\left\vert M\right\vert -\left\vert L\right\vert \leq cn\left\vert H^{\prime
}\left(  M\right)  \right\vert <\frac{c}{2}n\cdot n^{r-1},
\]
giving $\left\vert L\right\vert >\left(  c/2\right)  n^{r},$ as claimed.

Since $H^{\prime}\left(  L\right)  \subset H^{\prime}\left(  G\right)  $ and
\[
\left\vert H^{\prime}\left(  L\right)  \right\vert \geq\left\vert L\right\vert
/n>\left(  c/2\right)  n^{r}/n=\left(  c/2\right)  n^{r-1},
\]
the induction assumption implies that $H^{\prime}\left(  L\right)  $ covers an
induced subgraph $Z\subset G$ of type $H^{\prime}\left(  p,\ldots,p\right)  $
with $p=\left\lfloor c^{\left(  r-1\right)  ^{2}}\ln n\right\rfloor .$ Here we
use the inequalities
\[
n^{1-c^{r-2}}\geq n^{1-c}\geq n^{1/2}>2^{-4}\ln n\geq c^{\left(  r-1\right)
^{2}}\ln n.
\]

Since $H^{\prime}\left(  L\right)  $ covers $Z,$ there exist $R_{1}%
,\ldots,R_{p}\in H^{\prime}\left(  L\right)  $ such that $R_{1}\left(
H^{\prime}\right)  ,\ldots,R_{p}\left(  H^{\prime}\right)  $ are disjoint
subgraphs of $Z$ intersecting all its vertex classes. For every $i\in\left[
p\right]  ,$ let
\[
W_{i}=\left\{  v:\left(  \text{there exists }P\in L\text{ extending }%
R_{i}\right)  \text{ }\&\text{ }\left(  P\left(  v_{r}\right)  =v\right)
\right\}  .
\]

Write $d$ for the degree of $v_{r}$ in $H$ and note that each $v\in W_{i}$ is
joined to exactly $d$ vertices of $R_{i}\left(  H^{\prime}\right)  .$ Since,
by our selection, $d_{L}\left(  R_{i}\right)  \geq\left(  c/2\right)  n$ for
all $i\in\left[  p\right]  ,$ there is a set $X_{i}\subset W_{i}$ with
\[
\left\vert X_{i}\right\vert \geq\left(  cn/2\right)  /\binom{r-1}{d}\geq
cn/2^{r-1}%
\]
such that the vertices of $X_{i}$ have the same neighbors in $R_{i}\left(
H^{\prime}\right)  .$ Let $Y_{i}\subset\left[  r-1\right]  $ be the set of
classes of $Z$ containing the neighbors of the vertices of $X_{i}.$

Next, set $m=\left\lceil p/2^{r-2}\right\rceil ,$ and note that there is a set
$A\subset\left[  r-1\right]  $ with $\left\vert A\right\vert =m$ such that all
sets $Y_{i},$ $i\in A,$ are the same.

Define a bipartite graph $F$ with parts $A$ and $B=V\left(  G\right)  ,$
joining $i\in A$ to $v\in B$ if $v\in X_{i}.$ Since $\left\vert X_{i}%
\right\vert >\left(  c/2^{r-1}\right)  n$ for all $i\in A,$ we have%
\[
e\left(  F\right)  >\left(  c/2^{r-1}\right)  mn.
\]

Also, setting $s=\left\lfloor c^{r^{2}}\ln n\right\rfloor ,$ we find that%
\begin{align*}
s  &  \leq c^{r^{2}}\ln n=c^{2r-1}c^{\left(  r-1\right)  ^{2}}\ln n\leq\left(
c/2^{2r-2}\right)  \left\lfloor c^{\left(  r-1\right)  ^{2}}\ln n\right\rfloor
+1\\
&  \leq\left(  c/2^{r}\right)  p/2^{r-2}+1\leq\left(  c/2^{r}\right)  m+1.
\end{align*}

Therefore, by Lemma \ref{le1}, there exists $K_{2}\left(  s,t\right)  \subset
F$ with parts $S\subset A$ and $T\subset B$ such that $\left\vert S\right\vert
=s$ and $\left\vert T\right\vert =t>n^{1-c^{r-1}}.$

Let $G^{\prime}=G\left[  \cup_{i\in S}R_{i}\left(  H^{\prime}\right)  \right]
$ and $G^{\prime\prime}=G\left[  \cup_{i\in S}R_{i}\left(  H^{\prime}\right)
\cup T\right]  .$ Note that $G^{\prime}\subset Z,$ and clearly $G^{\prime}$ is
of type $H^{\prime}\left(  s,\ldots,s\right)  .$ Since each vertex $v\in T$ is
joined to exactly the same vertices of $\cup_{i\in S}R_{i}\left(  H^{\prime
}\right)  ,$ we see that $G^{\prime\prime}$ is of type $H\left(
s,\ldots,s,t\right)  .$

To finish the proof, we show that $L$ covers $G^{\prime\prime}.$ First, we see
that, for every edge $ij$ going across vertex classes of $G^{\prime\prime},$
there exists $h\in K_{2}\left(  L\right)  $ mapping some edge of $H$ onto
$ij$. Finally, taking $s$ distinct vertices $u_{1},\ldots,u_{s}\in T,$ by the
construction of $T$, for every $i\in S,$ there exists $P_{i}\in L$ with
$P_{i}|V\left(  X\right)  =R_{i}$ and $P_{i}\left(  v_{r}\right)  =u_{i}.$
Hence, $L$ covers $G^{\prime\prime},$ completing the induction step and the proof.
\end{proof}

\subsubsection*{Concluding remarks}

Using random graphs, it is easy to see that most graphs on $n$ vertices
contain substantially many copies of any fixed graph, but contain no
$K_{2}\left(  s,s\right)  $ for $s\gg\log n.$ Hence, Theorems \ref{th1},
\ref{th2}, and \ref{genZ} are essentially best possible.\bigskip

Finally, a word about the project mentioned in the introduction: in this
project we aim to give wide-range results that can be used further, adding
more integrity to extremal graph theory.

\end{document}